\documentclass[11pt, a4paper]{article}
\usepackage{amsmath,amssymb,latexsym,graphics,epsfig}
\usepackage{hyperref}
\usepackage{color}
\usepackage{amsthm}

\newtheorem{teorema}{Theorem}[section]
\newtheorem{lema}[teorema]{Lemma}
\newtheorem{corollari}[teorema]{Corollary}
\newtheorem{proposicio}[teorema]{Proposition}

\def\prova{{\boldmath  $Proof.$}\hskip 0.3truecm}

\def\final{\mbox{ \quad $\Box$}}

\def\A{\mbox{\boldmath $A$}}
\def\B{\mbox{\boldmath $B$}}

\def\dist{\mathop{\rm dist }\nolimits}

\def\excC{\mbox{$\varepsilon_{\Cpetita}$}}

\def\j{\mbox{\boldmath $j$}}

\def\A{\mbox{\boldmath $A$}}

\def\G{\Gamma}

\def\matrixI{\mbox{\boldmath $I$}}

\def\ev{\mbox{\rm ev}}
\def\sp{\mbox{\rm sp}}

\def\evG{\mbox{\rm ev\,$\Gamma$}}

\def\eve{\mbox{$\ev_e$\,$\Gamma$}}
\def\evE{\mbox{$\rm ev_{\Epetita}$\,$\Gamma$}}

\def\Cpetita{\mbox{{\tiny $C$}}}

\def\Epetita{\mbox{{\tiny $E$}}}

\def\ecc{\mathop{\rm ecc }\nolimits}

\def\tr{\mathop{\rm tr }\nolimits}

\def\ev{\mathop{\rm ev }\nolimits}
\def\sp{\mathop{\rm sp }\nolimits}

\begin{document}
\title{Edge-Distance-Regular Graphs\\ Are Distance-Regular
%\thanks{\noindent
%Research supported by the {\em Ministerio de Ciencia e Innovaci\'on}, Spain, and the {\em European Regional Development Fund} under project MTM2011-28800-C02-01 and by the {\em Catalan Research Council} under project 2009SGR1387.
%}
}

\author{M. C\'amara$^a$, C. Dalf\'{o}$^b$, C. Delorme$^c$, M.A. Fiol$^b$, H. Suzuki$^d$
\\ \\
{\small $^a$Tilburg University, Dept. of Econometrics and O.R.} \\
{\small  Tilburg, The Netherlands}\\
{\small (e-mail: {\tt
mcamara@ma4.upc.edu)}} \\
{\small $^b$Universitat Polit\`ecnica de Catalunya, BarcelonaTech} \\
{\small Dept. de Matem\`atica Aplicada IV, Barcelona, Catalonia}\\
{\small (e-mails: {\tt
\{cdalfo,fiol\}@ma4.upc.edu})} \\
{\small $^c$Universit\'e Paris Sud, LRI, Orsay CEDEX, France} \\
 {\small (e-mail: {\tt cd@lri.fr})}\\
{\small $^d$International Christian University,
Dept. of Mathematics and Computer Science} \\
{\small Mitaka, Tokyo 181-8585, Japan} \\
{\small (e-mail: {\tt hsuzuki@icu.ac.jp}) } }

%\date{}

\maketitle

\begin{abstract}
A graph is edge-distance-regular when it is distance-regular around
each of its edges and it has the same intersection numbers for any
edge taken as a root. In this paper we give some (combinatorial and algebraic) proofs of the fact
that every edge-distance-regular graph $\G$ is distance-regular and homogeneous.
More precisely, $\G$ is edge-distance-regular if and only if it is bipartite distance-regular or a
generalized odd graph. Also, we obtain the relationships between
some of their corresponding parameters, mainly, the distance polynomials and
the intersection numbers.
\end{abstract}

\section{Introduction}

In this paper, we use standard concepts and results about
distance-regular graphs (see, for example, Biggs \cite{b93}, or
Brouwer, Cohen, and Neumaier \cite{bcn89}), spectral graph theory
(see Cvetkovi\'c, Doob, and Sachs \cite{cds80}, or Godsil
\cite{g93}), and spectral and algebraic characterizations of
distance-regular graphs (see, for instance, Fiol \cite{f02}).

Edge-distance-regular graphs, introduced by Fiol and Garriga
\cite{fg01}, are analogous to dis\-tance-regular graphs but
considering the distance partitions induced by every edge instead of
each vertex. Thus, many known results for distance-regular graphs
have their counterpart for edge-distance-regular graphs such as, for
instance, the so-called spectral excess theorem. This theorem
characterizes (vertex- or edge-)distance-regular graphs by their
spectra and the (average) number of vertices at extremal distance (from every vertex or edge).
See Fiol and Garriga \cite{fg97} and C\'amara,  Dalf\'o, F\`abrega,
Fiol, and Garriga \cite{cdffg11} for the cases of distance-regular
and edge-distance-regular graphs, respectively. Also, for short
proofs, see Van Dam \cite{vd08} and Fiol, Gago, and Garriga
\cite{fgg09}.

A distance-regular graph $\G$ with diameter $d$ and odd-girth (that
is, the shortest cycle of odd length) $2d+1$ is called a {\em
generalized odd graph}, also known as an {\em almost-bipartite
distance-regular graph} or a {\em regular thin near $(2d+1)$-gon}.
The first name is due to the fact that odd graphs $O_k$ (see Biggs
\cite{b93}) are distance-regular and have such an odd-girth. Notice that, in this case, the
intersection parameters of $\G$ satisfy $a_0=\cdots=a_{d-1}=0$ and
$a_d\neq 0$. Recently, Van Dam and Haemers \cite{vdh10} showed that
any connected regular graph with $d+1$ distinct eigenvalues and
odd-girth $2d+1$ is a generalized odd graph. Moreover,
Lee and Weng \cite{lw11} used a variation of the spectral excess theorem for nonregular
graphs to show that, in fact, the regularity condition is not necessary, and Van Dam and Fiol \cite{vdf12}
gave a more direct short proof of the same result.

Here, we provide some (combinatorial and algebraic) proofs that, in
fact, any edge-distance-regular graph $\G$ is also distance-regular.
Moreover, if this is the case, $\G$ is either bipartite or a
generalized odd graph, and the relationship between the intersection
numbers of the corresponding distance partitions, induced by a
 vertex and by an edge, is made explicit.
Thus, a distance-regular graph $\G$ is edge-distance-regular if and
only if $\G$ is either bipartite or a generalized odd graph. In
fact, the `only if' part is also a consequence of a result by Martin
\cite{m92}, who proved that if a pair of vertices at distance $h$ is
a completely regular code in a distance-regular graph $\G$ with
diameter $d$, then either $h=1$ and $\G$ has intersection numbers
$a_1=\cdots = a_{d-1}=0$, or $h=d$ and $\G$ is antipodal.

In the rest of this section, we recall some concepts, terminology,
and results involved. Throughout this paper, $\G=(V,E)$ denotes a connected graph on $n=|V|$ vertices and $m=|E|$ edges, having adjacency matrix $\A$, and
spectrum $\sp \G=\{\lambda_0^{m_0},\ldots,\lambda_d^{m_d}\}$, where
$\lambda_0>\cdots>\lambda_d$, and the superscripts $m_i$ stand for the multiplicities.
The distance between two vertices $u,v$ is denoted by $\dist(u,v)$, so that the diameter of $\G$ is $D=\max\{\dist(u,v): u,v\in V\}$. Moreover, given $C\subset V$, the set $C_i=\G_i(C)=\{u\in V : \dist(u,C)=i\}$ is called the {\em $i$-th subconstituent} with respect to $C$,
where $\dist(u,C)=\min\{\dist(u,v): v\in C\}$ and $C_0=C$. In particular, when $C$ is a singleton, $C=\{u\}$, we write $\G_i(u)$ for $\G_i(\{u\})$ and set $\G(u)$ for $\G_1(u)$. The {\em eccentricity} or {\em covering radius} of $C$ is
$\ecc(C)=\excC=\max\{i : C_i\neq \emptyset\}$, so that we have the partition $V=C_0\cup C_1\cup\cdots\cup C_{\excC}$.

Given any two vertices $w,u$ at distance
$\dist(w,u)=i\ge 0$ of a graph $\G$,  we consider the numbers of
neighbors of $w$ at distance $i-1$, $i$, $i+1$  from $u$, that is,
$$
c_i(w,u)=|\G(w)\cap \G_{i-1}(u)|,\qquad a_i(w,u)=|\G(w)\cap
\G_i(u)|,\qquad b_i(w,u)=|\G(w)\cap \G_{i+1}(u)|,
$$
and $\G$ is distance-regular if these numbers only depend on $i$. In this case we write $c_i(w,u)=c_i$, $a_i(w,u)=a_i$ and  $b_i(w,u)=b_i$ and say that these numbers are {\em well defined}.
A useful characterization of distance-regularity is the existence of the so-called {\em distance-polynomials} $p_0,\ldots,p_d$ of $\G$ satisfying
\begin{equation}
\label{pi(A)=Ai}
p_i(\A)=\A_i,\qquad i=0,\ldots,d,
 \end{equation}
where $\A_i$ is the $i$-th distance matrix of $\G$, with entries $(\A_i)_{uv}=1$ if $\dist(u,v)=i$, and $(\A_i)_{uv}=0$ otherwise. Recall also that, if $\G$ is distance-regular,  the
intersection parameters $p_{ij}^k=|\G_i(u)\cap \G_j(v)|$, with
$\dist(u,v)=k$, for $i,j,k=0,\ldots,d$, are the Fourier
coefficients of the polynomial $p_ip_j$ in terms of the basis
constituted by the distance-polynomials of $\G$ with respect to the
scalar product
$$
\langle f,g\rangle = \frac{1}{n}\tr(f(\A)g(\A))
=\frac{1}{n}\sum_{i=0}^d m_i  f(\lambda_i)g(\lambda_i).
$$
 Thus, with $n_i=p_i(\lambda_0)=\|p_i\|^2$, we have the well-known relations
\begin{equation}\label{relating-int-param}
 n_k p_{ij}^k=\langle p_ip_j,p_k\rangle=\langle p_j,p_ip_k\rangle=n_j p_{ik}^j,
 \qquad i,j,k=0,\ldots,d.
\end{equation}
In particular, when $k=1$ and $i=j$, we have $n_1=\lambda_0=\delta$
(the degree of $\G$) and $p_{i1}^i=a_i$. Thus,
%\begin{equation}\label{relating-int-param-1}
% \delta p_{ii}^1=n_i a_i.
%\end{equation}
$\delta p_{ii}^1=n_i a_i$ and, hence, $p_{ii}^1=0$ if and only if $a_i=0$. In
fact notice that, for a general
graph, the condition $V_{i,i}(u,v)=|\G_i(u)\cap \G_i(v)|=\emptyset$ for any two adjacent vertices $u,v$ is equivalent to say that $a_i$ is well defined and null.

In an edge-distance-regular graph $\G=(V,E)$ with diameter $d$,
every pair of adjacent vertices $u,v\in V$ is a completely regular
code. More precisely, the distance partition
$\tilde{V_0},\ldots,\tilde{V}_{\tilde{d}}$ of $V$ induced by an edge
$uv\in E$, where $\tilde{V}_i=\tilde{V}_i(uv)=\G_i(uv)$ is the set of vertices at
distance $i$ from  $\{u,v\}$ and $\tilde{d}\in\{d-1,d\}$, is regular
and with the same edge-intersection numbers for any edge. That is, the numbers
$$
\tilde{a}_i(uv)=|\G(w)\cap \tilde{V}_i|,\qquad  \tilde{b}_i(uv)=|\G(w)\cap \tilde{V}_{i+1}|
,\qquad  \tilde{c}_i(uv)=|\G(w)\cap \tilde{V}_{i-1}|,
$$
do not depend neither on the edge $uv$ nor on the vertex $w\in\tilde{V}_i$, but only on the
distance  $i$, in which case we write them as $\tilde{a}_i,\tilde{b}_i,\tilde{c}_i$ for
$i=0,\ldots,\tilde{d}$  and say that they are well defined  (see C\'amara, Dalf\'o, F\`abrega, Fiol and Garriga \cite{cdffg11}
for more details).

\section{The characterization}
In the next result we show that every edge-distance-regular graph is
either bipartite distance-regular or a generalized odd graph.

%\begin{teorema}
%Let $\G$ be a distance-regular graph with diameter $D=d$ and intersection array
%$$
%\left(\begin{array}{c c c c c}
%0& c_1 & \cdots & c_{d-1}& c_d\\
% a_0& a_1 & \cdots & a_{d-1}& a_d\\
%  b_0&b_1 & \cdots & b_{d-1}& 0
%  \end{array}\right).
%  $$
%Then, $\G$ is edge-distance-regular if and only if it is either bipartite or a generalized odd graph.
%In this case, when $\G$ is a generalized odd graph, the edge-intersection array is
%\begin{equation}\label{relating-inter-arrays-nonbip}
%\left(\begin{array}{c c c c c}
%0& \tilde{c}_1 & \cdots &  \tilde{c}_{d-1}& \tilde{c}_d\\
%\tilde{a}_0& \tilde{a}_1 & \cdots &  \tilde{a}_{d-1}& \tilde{a}_d\\
%\tilde{b}_0&\tilde{b}_1 & \cdots &  \tilde{b}_{d-1}& 0\end{array}\right)
% =
%\left(\begin{array}{c c c c c }
%0& c_1 & \cdots &  c_{d-1}& 2c_d\\
%c_1& c_2-c_1 & \cdots &  c_d-c_{d-1}& a_d-c_d\\
%b_1& b_2& \cdots &  a_d & 0\end{array}\right),
%\end{equation}
%and, when $\G$ is bipartite:
%\begin{equation}\label{relating-inter-arrays-bip}
%\left(\begin{array}{c c c c c}
%0& \tilde{c}_1 & \cdots &  \tilde{c}_{d-2}& \tilde{c}_{d-1}\\
%\tilde{a}_0& \tilde{a}_1 & \cdots &  \tilde{a}_{d-2}& \tilde{a}_{d-1}\\
%\tilde{b}_0&\tilde{b}_1 & \cdots &  \tilde{b}_{d-2}& 0\end{array}\right)
% =
%\left(\begin{array}{c c c c c }
%0& c_1 & \cdots &  c_{d-2}& c_{d-1}\\
%c_1& c_2-c_1 & \cdots &  c_{d-1}-c_{d-2}& b_0-c_{d-1}\\
%b_1& b_2& \cdots &  b_{d-1} & 0\end{array}\right).
%\end{equation}
%\end{teorema}

\begin{teorema}\label{thm edr is dr}
Let $\G$ be a graph with diameter $d$. Then, the
following statements are equivalent:
\begin{itemize}
\item[$(a)$]
 $\G$ is edge-distance-regular; %with intersection
%parameters $(b_0,b_1,\ldots,b_{d-1}; c_1,c_2,\ldots,c_d)$..
 \item[$(b)$]
 $\G$ is distance-regular, either bipartite or a generalized odd graph.
\end{itemize}
Moreover, if this is the case and $\G$ has intersection array $\iota(\G)=\{b_0,b_1,\ldots, b_{d-1};$
$c_1,c_2,\ldots,c_d\}$ and it is bipartite, then its edge-intersection
array is
%\begin{equation}\label{relating-inter-arrays-bip}
%\left(\begin{array}{c c c c c}
%0& \tilde{c}_1 & \cdots &  \tilde{c}_{d-2}& \tilde{c}_{d-1}\\
%\tilde{a}_0& \tilde{a}_1 & \cdots &  \tilde{a}_{d-2}& \tilde{a}_{d-1}\\
%\tilde{b}_0&\tilde{b}_1 & \cdots &  \tilde{b}_{d-2}&
%0\end{array}\right)
% =
%\left(\begin{array}{c c c c c }
%0& c_1 & \cdots &  c_{d-2}& c_{d-1}\\
%c_1& c_2-c_1 & \cdots &  c_{d-1}-c_{d-2}& b_0-c_{d-1}\\
%b_1& b_2& \cdots &  b_{d-1} & 0\end{array}\right).
%\end{equation}
\begin{equation}\label{relating-inter-arrays-bip}
\tilde{\iota}(\G)=\{\tilde{b}_0,\tilde{b}_1,\ldots, \tilde{b}_{d-2};
\tilde{c}_1,\ldots,\tilde{c}_{d-2}, \tilde{c}_{d-1}\} = \{b_1, b_2,
\ldots  b_{d-1}; c_1, \ldots , c_{d-2}, c_{d-1}\},
\end{equation}
whereas, if $\G$ is a generalized odd graph, its edge-intersection
array is
\begin{equation}\label{relating-inter-arrays-nonbip}
\tilde{\iota}(\G)=\{\tilde{b}_0,\tilde{b}_1,\ldots, \tilde{b}_{d-1};
\tilde{c}_1,\ldots,\tilde{c}_{d-1}, \tilde{c}_{d}\} = \{b_1, b_2,
\ldots  a_{d}; c_1, \ldots , c_{d-1}, 2c_{d}\}.
\end{equation}

%\begin{equation}\label{relating-inter-arrays-nonbip}
%\left(\begin{array}{c c c c c}
%0& \tilde{c}_1 & \cdots &  \tilde{c}_{d-1}& \tilde{c}_d\\
%\tilde{a}_0& \tilde{a}_1 & \cdots &  \tilde{a}_{d-1}& \tilde{a}_d\\
%\tilde{b}_0&\tilde{b}_1 & \cdots &  \tilde{b}_{d-1}& 0\end{array}\right)
% =
%\left(\begin{array}{c c c c c }
%0& c_1 & \cdots &  c_{d-1}& 2c_d\\
%c_1& c_2-c_1 & \cdots &  c_d-c_{d-1}& a_d-c_d\\
%b_1& b_2& \cdots &  a_d & 0\end{array}\right),
%\end{equation}
\end{teorema}

\begin{figure}[h]
\begin{center}
\includegraphics[width=15cm]{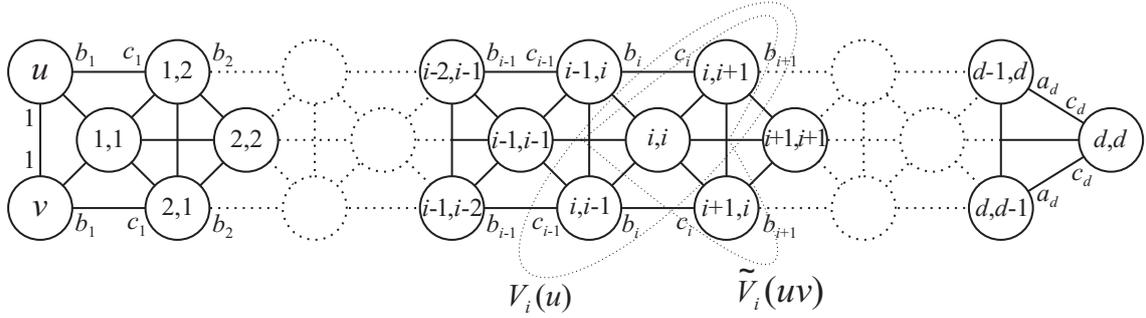}
\vskip -7cm \caption{The distance partition induced by two adjacent
vertices. (Every set $V_{i,j}$ is represented by its subindices.) }
\label{edge-diagram}
\end{center}
\end{figure}

\prova
As the complete graphs clearly satisfy the result, we can
assume that $d>1$. Given two adjacent vertices $u$ and $v$ of $\G$, let us
consider the intersection numbers $ p_{ij}^k(u,v)=|\G_i(u)\cap
\G_j(v)|$, so that the vertex partition induced by the distances
from $u$ and $v$ is shown in Fig. \ref{edge-diagram}, where
$V_{i,j}=V_{i,j}(u,v)=\G_i(u)\cap \G_j(v)$. (Notice that $V_{i,j}=\emptyset$ when
$|i-j|>1$, as $\dist(u,v)=1$.) Let
$\tilde{V}_{0},\ldots,\tilde{V}_{\tilde{d}}$ be the distance
partition induced by the edge $uv$, and define $\tilde{a}_i(uv)$,
$\tilde{b}_i(uv)$ and $\tilde{c}_i(uv)$ as above. Clearly,
$\tilde{a}_0(uv)=\tilde{a}_0=c_1=1$.

%

%\noindent{\bf Some relations.}\\
We have the following facts:

\begin{itemize}
\item[$(i)$]
For $i=1,\ldots,d-1$, all vertices adjacent to $w\in
V_{i,i-1}\subset \tilde{V}_{i-1}$ and at distance $i+1$ from $u$
are in $V_{i+1,i}\subset \tilde{V}_{i}$. Hence,
\begin{eqnarray}
|\G(w)\cap \tilde{V}_{i}| & = & |\G(w)\cap \G_{i+1}(u)|+|\G(w)\cap V_{i,i}|\nonumber\\
  & = & b_i(w,u)+|\G(w)\cap V_{i,i}|,\qquad i=1,\ldots,d-1.\label{desig1}
\end{eqnarray}
%Analogously, if $w\in V_{i-1,i}\subset \tilde{V}_{i-1}$, we get
%\begin{eqnarray}\label{desig2}
%|\G(w)\cap \tilde{V}_{i}| & = & b_i(w,v)+|\G(y)\cap V_{i,i}|,\qquad
%i=1,\ldots,d-1.
%\end{eqnarray}

\item[$(ii)$]
For $i=1,\ldots,d-1$, all vertices adjacent to $w\in
V_{i,i+1}\subset \tilde{V}_{i}$ and at distance $i-1$ from $u$ are in $V_{i-1,i}\subset \tilde{V}_{i-1}$. Thus,
\begin{equation}
\label{desig3} |\G(w)\cap \tilde{V}_{i-1}| =|\G(w)\cap \G_{i-1}(u)|
= c_i(w,u),\qquad i=1,\ldots,d-1.
\end{equation}
Moreover, assuming that $a_i$ is well defined with value $a_i=0$, the vertices adjacent to vertex $w$ but at distance
$i$ from $v$ are in $V_{i-1,i}\cup V_{i+1,i}$ since $V_{i,i}=\emptyset$, whereas $|\G(w)\cap V_{i,i+1}|=0$ since $w\in V_{i,i+1}\subset \G_{i}(u)$. Consequently,
\begin{eqnarray}
\hskip -.5cm |\G(w)\cap \tilde{V}_{i-1}|+|\G(w)\cap
\tilde{V}_{i}| & = & |\G(w)\cap (V_{i-1,i}\cup
V_{i+1,i})| \nonumber \\
 & = & |\G(w)\cap \G_{i}(v)| = c_{i+1}(w,v),\qquad
i=1,\ldots,d-1. \label{desig3b}
\end{eqnarray}
%Similarly, if $w\in V_{i+1,i}\subset \tilde{V}_{i}$, we have
%\begin{equation}
%\label{desig4} |\G(w)\cap \tilde{V}_{i-1}|=c_i(w,v),\qquad
%i=1,\ldots,d-1.
%\end{equation}

\item[$(iii)$]
For $i=1,\ldots,d$, all vertices adjacent to $w\in
V_{i,i}\subset\tilde{V}_{i}$ and at distance $i-1$ from $u$ are
either in $V_{i-1,i}$ or $V_{i-1,i-1}$. Thus,
$$
|\G(w)\cap V_{i-1,i}|+|\G(w)\cap V_{i-1,i-1}|=|\G(w)\cap
\G_{i-1}(u)|=c_i(w,u),
$$
and, analogously,
$$
|\G(w)\cap V_{i,i-1}|+|\G(w)\cap V_{i-1,i-1}|=|\G(w)\cap
\G_{i-1}(w)|=c_i(w,v).
$$
Therefore,
\begin{eqnarray}
|\G(w)\cap \tilde{V}_{i-1}| & = &
 c_i(w,u)+c_i(w,v)-|\G(w)\cap V_{i-1,i-1}|,\qquad i=1,\ldots,d.\label{desig5}
\end{eqnarray}
\end{itemize}

\noindent{$(b)\Rightarrow (a)$:}\\
Assume that $\G$ is distance-regular with intersection parameters $a_i,b_i,c_i$, for $i=0,\ldots,d$
(the parameters $b_i$ and $c_i$ are indicated in Fig. \ref{edge-diagram}).
If $\G$ is either bipartite or a generalized odd graph, then $a_i=0$
and $V_{i,i}=\emptyset$ for every $i=1,\ldots,d-1$ (in particular, as $a_{d-1}=0$, the parameters $a_d$ and $c_d$ are those indicated in the same figure for the nonbipartite case; otherwise, we also have $V_{d,d}=\emptyset$ and, hence, $a_d=c_d=0$). Thus, from
the above reasonings we obtain:

\begin{itemize}
\item[$(i)$] By Eq. (\ref{desig1}) %and (\ref{desig2})
and for every vertex
$w\in \tilde{V}_{i-1}$, we have
$$
\tilde{b}_{i-1}(uv)=|\G(w)\cap \tilde{V}_{i}|=b_i.
$$
Hence,   $\tilde{b}_{i}$ is well defined for $i=0,\ldots,d-2$.
Moreover,
if $\G$ is a generalized odd graph, $\tilde{b}_{d-1}=a_d$.
%and $\tilde{b}_d=0$,
%whereas, if $\G$ is bipartite, $\tilde{b}_{d-1}=0$.

\item[$(ii)$]
 By Eq. (\ref{desig3}) %and (\ref{desig4})
and for every vertex $w\in \tilde{V}_{i}$,
$$
\tilde{c}_i(uv)=|\G(w)\cap \tilde{V}_{i-1}|=c_i.
$$
Consequently, $\tilde{c}_i$ is well defined for $i=1,\ldots,d-1$.

\item[$(iii)$]
Moreover, if $\G$ is a generalized odd graph,
$\tilde{V}_{d,d}\neq\emptyset$, and Eq. (\ref{desig5}) yields that, for
every vertex $w\in \tilde{V}_{d,d}$,
$$
\tilde{c}_d(uv) = |\G(w)\cap \tilde{V}_{d-1}|=2c_d.
$$
%%all vertices adjacent to $w\in V_{dd}=\tilde{V}_d$ and at distance $d-1$ from $u$ (respectively,
% from $v$) are at  $V_{d-1,d}$ (respectively, at $V_{d,d-1}$), as $V_{d-1,d-1}=\emptyset$.
%  Hence, $|\G(w)\cap V_{d-1,d}|=|\G(w)\cap V_{d,d-1}|=c_d$ and $\tilde{c}_d=2 c_d$.
%%Moreover,
%%if $\G$ is a generalized odd graph, $\tilde{c}_d=2 c_d$ since, for every $w\in V_{dd}$,
%$|\G(w)\cap V_{d-1,d}|=|\G(w)\cap V_{d,d-1}|=c_d$.
\end{itemize}
Summarizing, all intersection numbers $\tilde{b}_i$, $i=0,\ldots,\tilde{d}-1$,
and $\tilde{c}_i$, $i=1,\ldots,\tilde{d}$ are well defined, and $\G$ is edge-distance-regular.
(Of course, the other intersection numbers
are just $\tilde{a}_i=\delta-\tilde{b}_i-\tilde{c}_i$, for $i=0,\ldots,\tilde{d}$.)
% Finally, its edge-intersection array is as in (\ref{relating-inter-arrays-nonbip}) and
% (\ref{relating-inter-arrays-bip})
%since the equation $\tilde{a}_i=\delta-\tilde{b}_i-\tilde{c}_i$ yields
%$\tilde{a}_0=\delta-b_1=c_1$,
%$\tilde{a}_i=\delta-c_i-b_{i+1}=c_{i+1}-c_i$ for $i=1,\ldots,d-1$,
%$\tilde{c}_d=\delta-2c_d=a_d-c_d$ (for the generalized odd graphs);
%and $\tilde{a}_i=c_{i+1}-c_i$ for $i=1,\ldots,d-2$, $\tilde{a}_{d-1}=\delta-c_{d-1}$
%(when $\G$ is bipartite).

\noindent{$(a)\Rightarrow (b)$:}\\
Assume that $\G$ is
edge-distance-regular with edge-intersection numbers
$\tilde{a}_i,\tilde{b}_i,\tilde{c}_i$, for $i=0,\ldots,\tilde{d}$. Then,
to show that the numbers $c_i(w,u)$, $a_i(w,u)$ and $b_i(w,u)$
depend only on the distance $i=\dist(w,u)$, we proceed by induction on
$i$. To begin with, observe that $b_0$, $a_0$, and $c_1$ are well defined, with values $b_0=b_0(u,u)=\tilde{b}_0+\tilde{a}_0=\tilde{b}_0+1$, $a_0=a_0(u,u)=0$,
and $c_1=c_1(w,u)=1$.

Now, assume that, for some  $i\le d-1$, $c_i$ and $a_{i-1}$ are well defined with
$a_{i-1}=0$ or, equivalently,  $V_{i-1,i-1}=\emptyset$.
Thus, in order to show that $c_{i+1}$  and $a_i$ exist,
consider a shortest path $v,u,\ldots,w$ of length $i+1(\le d)$, so that
$w\in V_{i,i+1}(u,v)$. Then, Eq. (\ref{desig3}) gives
$$
\tilde{c}_i=|\G(w)\cap\tilde{V}_{i-1}|=c_i.
$$
 Now, suppose that $u'v'$ is an arbitrary edge. If we assume that $V_{i,i}(u',v')\neq \emptyset$, there exists a vertex $w'\in V_{i,i}(u',v')$ and Eq. (\ref{desig5}) gives
 $$
 \tilde{c}_i=|\G(w)\cap
\tilde{V}_{i-1}|=2c_{i},
$$
which would be a contradiction. Thus $V_{i,i}(u',v')=\emptyset$ and
$a_i$ exists, being $a_i=0$. Moreover, coming back to the edge $uv$,  Eqs. (\ref{desig3}) and (\ref{desig3b}) yield
$$
c_{i+1}(w,v)=|\G(w)\cap \G_{i-1}(u)| +|\G(w)\cap \tilde{V}_{i}|=c_i+\tilde{a}_i,
$$
and, since $w$ and $v$ are arbitrary vertices at distance $i+1$, $c_{i+1}$ is well defined.

Thus, the numbers
$b_i=b_0-c_i-a_i$, for
$i=1,\ldots,d-1$, are also well defined; and the same holds for $a_d=b_0-c_d$ when $\G$ is nonbipartite. Consequently, $\G$ is as claimed.
\final

\subsection{Homogeneous graphs}
As a consequence of the previous results, we next prove that the edge-distance-regular graphs are, in fact, a particular case of the so-called (1-)homogeneous graphs introduced by Nomura in \cite{n94}.
A graph $\G$ is called {\em homogeneous} if for all nonnegative integers $r,s,i,j$ and pairs of edges $uv$, $u'v'$,
$$
x\in V_{s,r}(u,v),\ x'\in V_{s,r}(u',v')\Longrightarrow |\G(x)\cap V_{i,j}(u,v)|=|\G(x')\cap V_{i,j}(u',v')|,
$$
where $V_{s,r}(u,v)$ is defined as before. In other words, $\G$ is homogeneous if the partition of Fig.~\ref{edge-diagram} is a regular (or equitable) partition with the same intersection numbers for every pair of adjacent vertices $u,v$. For a detailed study of regular partitions, see Godsil \cite{g93} and Godsil and McKay \cite{gmk80}.

In  \cite[Lemma~5.3]{cn05}, Curtin and Nomura showed that every bipartite or almost bipartite distance-regular graph is homogeneous. Then, from Theorem \ref{thm edr is dr} we have:

\begin{corollari}
 Every edge-distance-regular graph  is homogeneous. \final
\end{corollari}

\begin{figure}[h]
\begin{center}
\includegraphics[width=16cm]{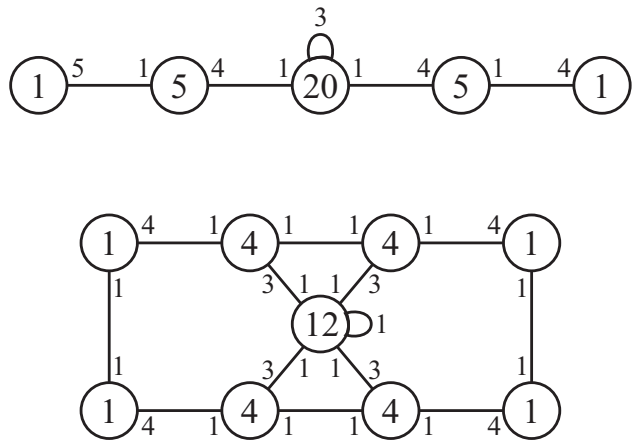}
\vskip -14cm
\caption{The intersection diagrams of the Wells graph as a distance-regular and homogeneous graph.}
\label{wells}
\end{center}
\end{figure}

However, the converse is not true. A counterexample is, for instance, the Wells graph $W$, which is the unique distance-regular graph on $n=32$ vertices and intersection array $\iota(W)=\{5,4,1,1;1,1,4,5\}$ (see Brouwer, Cohen and Neumaier \cite[Theorem~9.2.9]{bcn89}). Thus, since $a_2=3$ and $a_i=0$ for every $i\neq 2$, $W$ is homogeneous (see Nomura \cite[Lemma~2]{n94}) but not edge-distance-regular. The intersection diagrams induced by a vertex and two adjacent vertices, showing that $W$ is homogeneous, are illustrated in Fig.\ref{wells}.

%

%%%%%%%%%%%%%%%%%%%%%%%%%%%%%%%%%%%%%%%%%%%%%%%%
\section{An algebraic approach}
%%%%%%%%%%%%%%%%%%%%%%%%%%%%%%%%%%%%%%%%%%%%%%%%
%%%%%%%%%%%%%%%%%%%%%%%%%%%%%%%%%%%%%%%%%%%%%%%%
The matrix approach to the result of Theorem \ref{thm edr is dr} is based on the relationship of orthogonal polynomials with (edge-)distance-regularity of graphs. Recall that every system of orthogonal polynomials of a discrete variable $\{r_i\}_{0\leq i\leq d}$ satisfies a three term recurrence of the form
\begin{equation}\label{eqn three term rec}
xr_i=\beta_{i-1}r_{i-1}+\alpha_i r_i+\gamma_{i+1}r_{i+1},\qquad 0\leq i\leq d,
\end{equation}
where the terms $b_{-1}$ and $c_{d+1}$ are zero.

In \cite{fgy96b} Fiol, Garriga and Yebra gave the following characterization of distance-regularity.
\begin{teorema}
\label{teo caract d-r}
A  regular graph  with  $d+1$ distinct eigenvalues is distance-regular if
its distance matrix $\A_d$ is a polynomial of degree $d$ in $\A$, that is,
\begin{equation}
\label{pd(A)=Ad}
p_{d}(\A)=\A_{d}.
\end{equation}
\end{teorema}

Moreover, such a polynomial is the highest degree member of the predistance polynomials  $\{p_i\}_{0\leq i\leq d}$ of $\Gamma$---which are orthogonal with respect to the scalar product $\langle p,q\rangle=\frac{1}{n}\tr(p(\A)q(\A))$, and are normalized in such a way that $\|p_i\|^2=p_i(\lambda_0)$. When the graph is distance-regular, these polynomials satisfy a three term recurrence as in Eq. (\ref{eqn three term rec}) with $\{\gamma_i,\alpha_i,\beta_i\}=\{c_i,a_i,b_i\}$.

The first ingredient for the algebraic approach to our results is the following proposition from C\'amara, Dalf\'o, F\`abrega, Fiol, and Garriga \cite{cdffg11}, which describes the different possibilities for the edge-diameter and the local spectra of edges in an edge-distance-regular graph. In a regular graph, the $e$-local spectrum of an edge $e=uv$ is constituted by those eigenvalues $\lambda_i\in \eve$ such that the orthogonal projection of the characteristic vector of $\{u,v\}$ on their corresponding $\lambda_i$-eigenspace is not null; see  \cite{cdffg11} for more details. The set of these {\em $e$-local eigenvalues}
is denoted by $\eve$.

\begin{proposicio}
\label{prop e-d-r propietats}
Let $\G$ be an edge-distance-regular graph with diameter $D$ and $d+1$ distinct eigenvalues. Then, the following statements hold:
\begin{itemize}
\item [$(a)$] $\G$ is regular.
\item [$(b)$] $\G$ has spectrally maximum diameter $(D=d\/)$ and its edge-diameter satisfies:
  \begin{itemize}
  \item [$(b1)$] If $\G$ is nonbipartite, then $\tilde{D}=d$;
  \item [$(b2)$] If $\G$ is bipartite, then $\tilde{D}=d-1$.
  \end{itemize}
\item [$(c)$] $\G$ is edge-spectrum regular and, for every
$e\in E$, the $e$-spectrum satisfies:
\begin{itemize}
\item [$(c1)$] If $\G$ is nonbipartite, then $\eve=\evG$;
\item [$(c2)$] If $\G$ is bipartite, then $\eve=
\evG\setminus \{-\lambda_0\}$ .
\end{itemize}
\end{itemize}
\end{proposicio}

Let $\ev_{\Epetita}\G=\bigcup_{e\in E}\eve$,  $\ev^*_{\Epetita}\G=\evE\setminus \{-\lambda_0\}$ and $\tilde{d}=|\ev^*_{\Epetita}\G|$.  Then, by using this notation, Proposition \ref{prop e-d-r propietats} establishes that, if $\G$ is edge-distance-regular and nonbipartite, then $\evE=\eve=\evG$ for every edge $e\in E$, and $\evE=\eve=\evG\setminus\{-\lambda_0\}$ otherwise.
In both cases all edges have the same local spectrum ($\G$ is said to be edge-spectrum regular), with $\tilde{d}+1=|\evE|$ distinct $e$-local eigenvalues.

The role of the distance matrices in the study of edge distance-regularity is played by the distance incidence matrices. More precisely, for a given graph $\G=(V,E)$, let $\tilde{D}=\max_{uv\in E}\{\ecc(\{u,v\})\}$ be its {\em edge-diameter}. Then, for every $i=0,1,\ldots,\tilde{D}$, the $i$-{\it incidence matrix} of $\G$ is the $(|V|\times
|E|)$-matrix $\B_i=(b_{ue})$ with entries $b_{ue}=1$ if
$\dist(u,e)=i$, and $b_{ue}=0$ otherwise. The following result corresponds to Theorems 4.4 and 4.9 of \cite{cdffg11}.

\begin{teorema}
\label{teo caract e-d-r}
A regular graph $\G=(V,E)$ with edge-diameter $\tilde{D}$ and $\tilde{d}+1=|\evE|$
is edge-distance-regular if and only if any of the following conditions hold:
\begin{itemize}
\item[$(a)$]
There exists a polynomial $\tilde{p}_i$ of
degree $i$ such that
\begin{equation}
\label{tilde-pi(A)}
\tilde{p}_i(\A)\B_0=\B_i, \qquad i=0,1,\ldots, \tilde{D}.
\end{equation}

\item[$(b)$]
$|\G_{\tilde{d}}(e)|=2\tilde{p}_{\tilde{d}}(\lambda_0)$ for every edge $e\in E$.
 \end{itemize}
 \end{teorema}
Moreover, if this is the case,  these  polynomials are the edge-predistance-polynomials, $\{\tilde{p}_i\}_{0\leq i\leq \tilde{d}}$, with $\tilde{d}=\tilde{D}$, satisfying a three term recurrence as in Eq. (\ref{eqn three term rec}) with $\{\gamma_i,\alpha_i,\beta_i\}=\{\tilde{c}_i,\tilde{a}_i,\tilde{b}_i\}$. Here we use the following consequence, which can be seen as the analogue of Theorem \ref{teo caract d-r} for edge-distance-regularity.

\begin{corollari}
\label{coro-edr}
A regular graph $\G=(V,E)$ with $\tilde{d}+1=|\evE|$
is edge-distance-regular if and only if
\begin{equation}
\label{edge-pd(A)=Ad}
\tilde{p}_{\tilde{d}}(\A)\B_0=\B_{\tilde{d}}.
\end{equation}
\end{corollari}
\prova
We only need to prove sufficiency, which is straightforward if we multiply both sides of Eq. (\ref{edge-pd(A)=Ad})
by the all-1 vector $\j$ and apply Theorem \ref{teo caract e-d-r}$(b)$.
\final

\subsection{A matrix approach to Theorem~\ref{thm edr is dr}}
Now, we can rewrite Theorem~\ref{thm edr is dr} and prove it by showing that the distance polynomials exist if and only if the edge-distance-polynomials do.

\begin{teorema}\label{thm edr is dr algebraic}
Let $\G$ be a graph with diameter $d$ and edge-diameter $\tilde{d}\in\{d-1,d\}$. Then, the
following statements are equivalent:
\begin{itemize}
\item[$(a)$]
$\G$ is distance-regular with $a_0=\cdots=a_{d-1}=0$; %, either bipartite or a generalized odd graph.
\item[$(b)$]
$\G$ is edge-distance-regular. %with intersection
%parameters $(b_0,b_1,\ldots,b_{d-1}; c_1,c_2,\ldots,c_d)$..
\end{itemize}
Moreover, if this is the case, the relationships between the corresponding edge-distance-polynomials and distance polynomials are:
\begin{eqnarray}
\tilde{p}_i & = & p_i-p_{i-1}+p_{i-2}-\cdots + (-1)^ip_0,\label{p-tilde-i-vs-p's} \\
p_i & = & \tilde{p}_i + \tilde{p}_{i-1},\label{p-i-vs-p-tildes's}
\end{eqnarray}
for $i=0,\ldots,d-1$, and, if $\G$ is nonbipartite,
\begin{eqnarray}
\tilde{p}_d & = & \textstyle \frac{1}{2}(p_d-p_{d-1}+p_{d-2}-\cdots + (-1)^dp_0);\label{p-tilde-d-vs-p's} \\
p_d & = & 2\tilde{p}_d + \tilde{p}_{d-1},\label{p-d-vs-p-tildes's}
\end{eqnarray}
whereas, if $\G$ is bipartite,
\begin{eqnarray}
p_d & = & H-\tilde{q}_{d-1}-\tilde{q}_{d-2}, \label{bip-p-d-vs-p-tildes's}
\end{eqnarray}
where $H=p_0+\cdots+p_d$ is the Hoffman polynomial satisfying $H(\lambda_i)=\delta_{ij}n$,
and $\tilde{q}_{i}=\tilde{p}_{0}+\cdots+\tilde{p}_{i}$, for $i=d-2,d-1$.
\end{teorema}

\prova
\noindent $(a)\Rightarrow (b)$: \\
Since $a_i=0$  we have that
\begin{equation}
\label{AiB0}
\A_i\B_0=\B_i+\B_{i-1}, \qquad i=0,\ldots,d-1.
\end{equation}
Moreover, when $\G$ is nonbipartite, $a_d\neq 0$ and
$$
\A_d\B_0=2\B_d+\B_{d-1}.
$$
Then, multiplying both sides of Eq. (\ref{pd(A)=Ad}) by $\B_0$ (on the right) and using all the above equations, we get
\begin{eqnarray*}
p_d(\A)\B_0 & = & 2\B_d+\B_{d-1}\\
 & = & 2\B_d+\A_{d-1}\B_0-\B_{d-2} \\
  & = &  2\B_d+\A_{d-1}\B_0-\A_{d-2}\B_0+\B_{d-3} \\
  & \vdots &  \\
  & = & 2\B_d+(\A_{d-1}-\A_{d-2}+\A_{d-3}-\cdots +(-1)^{d+1}\A_0)\B_0.
\end{eqnarray*}
Thus, using Eq. (\ref{pi(A)=Ai}), we get that
$\tilde{p}_d(\A)\B_0=\B_d$ with
$$
\textstyle
\tilde{p}_d =  \frac{1}{2}(p_d-p_{d-1}+p_{d-2}-\cdots + (-1)^dp_0),
 $$
as claimed in Eq. (\ref{p-tilde-d-vs-p's}),
and $\G$ is edge-distance-regular by Corollary \ref{coro-edr}.

Similarly, by using again Eq. (\ref{AiB0}), we get
\begin{eqnarray*}
 p_{i}(\A)\B_0 & = & \B_{i}+\B_{i-1} \\
     & = & \B_{i}+\A_{i-1}\B_0-\B_{i-2} \\
    & \vdots &  \\
    & = & \B_{i}+(\A_{i-1}-\A_{i-2}+\cdots +(-1)^{i}\A_0)\B_0,
 \end{eqnarray*}
so that
$$
\tilde{p}_{i}=p_{i}-p_{i-1}+\cdots + (-1)^{i}p_0,\qquad i=0,\ldots,d-1,
 $$
as claimed in Eq. (\ref{p-tilde-i-vs-p's}).
Moreover, if $\G$ is bipartite, then $a_d=0$, $\tilde{d}=d-1$, and
the existence of the edge-distance-polynomial $\tilde{p}_{d-1}$ imply that $\G$ is edge-distance-regular, again by Corollary \ref{coro-edr}.

Eqs. (\ref{p-i-vs-p-tildes's}) and (\ref{p-d-vs-p-tildes's}) are immediate consequences from Eqs. (\ref{p-tilde-i-vs-p's}) and (\ref{p-tilde-d-vs-p's}). Finally, when $\G$ is bipartite, we can obtain $p_d$ by adding all the equalities in Eq. (\ref{p-i-vs-p-tildes's}) to get $q_{d-1}=p_0+\cdots+p_{d-1}=\tilde{q}_{d-1}+\tilde{q}_{d-2}$, so that
$p_d=H-q_{d-1}$, where $H$ is the Hoffman polynomial (see Hoffman \cite{hof63}), and we have Eq. (\ref{bip-p-d-vs-p-tildes's}).

\noindent{$(b)\Rightarrow (a)$:}\\
%Since $\G$ is edge-distance-regular, we have $\tilde{p}_{\tilde{d}}(\A)\B_0=\B_{\tilde{d}}$.
Note first that $\B_0\B_0^{\top}=\A+\delta\matrixI$, where $\delta$ is the degree of $\Gamma$. Now, we show that all distance-polynomials exist, in particular $p_d(\A)=\A_d$, and apply Theorem~\ref{teo caract d-r}.

Suppose that $\G$ is nonbipartite. By Proposition~\ref{prop e-d-r propietats} we know that $\eve=\evG$ and, consequently, $\tilde{d}=d$. As in Theorem \ref{thm edr is dr}, the proof is by induction on $i$.
First, $a_1=0$ (as $d>1$), $c_1=1$, and the first two distance polynomials exist $p_0=1$ and $p_1=x$.
Now, assume that, for some  $i< d-1$,  $a_i=0$, $c_i$ is well defined, and there exist the polynomials $p_{i-1}$ and $p_i$.

To compute the product $\B_{i}\B_0^{\top}$, let us consider two vertices $u,w$ at distance $\dist(u,w)=i<d-1$
and take $v\in \G(u)\cap \G_{i-1}(w)$. Then, $(\B_{i}\B_0^{\top})_{uw}=\sum_{e\in E}(\B_{i})_{ue}(\B_0^{\top})_{ew}$. Each term of the sum equals $1$ for every vertex $w'\in \G(w)$ such that the edge $e=ww'$ is at distance $i$ from $u$ (since $(\B_i)_{u,ww'}=(\B_0^{\top})_{ww',w}=1$); that is, for every vertex
$$
w'\in \G(w)\cap V_{i,i-1}(u,v)\cup V_{i,i}(u,v)\cup V_{i+1,i}(u,v)=\G(w)\cap V_{i+1,i}(u,v),
 $$
where we used that, since $a_i=0$, $\G(w)\cap V_{i,i-1}(u,v)=\emptyset$ (notice that if $w'\in V_{i,i-1}(u,v)$, then we would have $a_i(w,u)\neq 0$) and $V_{i,i}(u,v)=\emptyset$.
Hence,  $(\B_{i}\B_0^{\top})_{uw}=\tilde{b}_{i-1}$.
Reasoning similarly, for a vertex $w$ such that  $\dist(u,w)=i+1$, and $v\in \G(u)\cap \G_{i}(w)$ (that is, $w\in V_{i+1,i}(u,v)\subset \tilde{V_i}$), we have $(\B_{i}\B_0^{\top})_{uw}=\tilde{a}_i+c_{i}$.
Otherwise, if $\dist(u,w)\neq i,i+1$, then $(\B_{i}\B_0^{\top})_{uw}=0$. Consequently,
$$
\B_{i}\B_0^{\top}=\tilde{b}_{i-1}\A_i+(\tilde{a}_i+c_{i})\A_{i+1},
$$
and, multiplying Eq. (\ref{tilde-pi(A)}) by $\B_0^{\top}$ on the right, we get
$$
\tilde{p}_i(\A)(\A+\delta\matrixI)=\tilde{b}_{i-1}p_i(\A)+(\tilde{a}_i+c_{i})\A_{i+1}.
$$
Thus, the distance polynomial of degree $i+1$ is
\begin{equation}
\label{p_(i+1)}
\textstyle
p_{i+1}=\frac{1}{\tilde{a}_i+c_i}((x+\delta)\tilde{p}_{i}-\tilde{b}_{i-1}p_i),
\end{equation}
graph $\G$ is {\em $(i+1)$-partially distance-regular} (that is, $p_j(\A)=\A_j$ for $j=0,\ldots,i+1$), and $c_{i+1}$ is well defined; see Dalf\'o, Van Dam, Fiol, Garriga, and Gorissen \cite{ddfgg11}.
(Notice that this could also be derived from Eqs. (\ref{desig3}) and (\ref{desig3b}) yielding $c_{i+1}(w,u)=|\G(w)\cap \G_{i-1}(v)| +|\G(w)\cap
\tilde{V}_{i}|=c_i+\tilde{a}_i$, as in the proof of Theorem \ref{thm edr is dr}.)
Now, let $u,v$ be two arbitrary adjacent vertices.
If $w\in V_{i+1,i+1}\neq \emptyset$, Eq.
(\ref{desig5}) yields  $\tilde{c}_{i+1}=|\G(w)\cap
\tilde{V}_{i}|=2c_{i+1}> c_{i+1}$.
It follows that if a vertex is at distance at most $i+1$ from an end of an edge, it is at distance at most $i$ from the other end. Thus, the diameter of $\G$ is $i+1<d$, which is a contradiction.
Hence, $V_{i+1,i+1}(u,v)=\emptyset$ and,
since $u,v$ are generic vertices,
%from (\ref{relating-int-param-1}), we conclude that
$a_{i+1}=0$ is well defined.

The induction above proves that there exist all the distance-polynomials $p_0,p_1,\ldots,p_{d-1}$ and, also,
we have $p_d=H-p_0-p_{1}-\cdots-p_{d-1}$. Then, $\G$ is distance-regular with $a_0=\cdots=a_{d-1}=0$.

%Now, to compute $\A\tilde{p}_{d}(\A)$, we use the three term recurrence in Eq. (\ref{eqn three term rec}), satisfied by the edge-predistance polynomials $\{\tilde{p}_k\}_{0\leq k\leq d}$,
%$$
%\tilde{b}_{d-1}\tilde{p}_{d-1}(\A)+\tilde{a}_d\tilde{p}_{d}(\A)
%+\delta\tilde{p}_{d}(\A)=\tilde{b}_{d-1}\A_d.
%$$
%Then, the polynomial $p=(\tilde{a}_d
%+\delta)(1/\tilde{b}_{d-1})\tilde{p}_{d}+\tilde{p}_{d-1}$ satisfies $p(\A)=\A_d$ and $\G$ is distance-regular
%by Theorem \ref{teo caract d-r}.
%Besides, using the relationships between the intersection numbers given in Theorem \ref{thm edr is dr} and that $\delta=a_d+c_d$, the above coefficient of $\tilde{p}_{d}$ equals $2$ and $p_d=p=2\tilde{p}_{d}+\tilde{p}_{d-1}$,
%according to (\ref{p-d-vs-p-tildes's}).

%The other equalities in (\ref{p-i-vs-p-tildes's}) giving the polynomials $p_i$ are proved analogously by using
%that
%$$
%\B_{i}\B_0^{\top}=\tilde{b}_{i-1}\A_i+\tilde{c}_{i}\A_{i+1},\qquad i=0,\ldots,d-1.
%$$

Suppose now that $\G$ is bipartite. Then $\tilde{d}=d-1$ and $\eve=
\evG\setminus \{-\lambda_0\}$. Now, reasoning similarly as above from a shortest path $u,v,\ldots,w$ of length $d$, we have
$$
\B_{d-1}\B_0^{\top}=\tilde{b}_{d-2}\A_{d-1}+(\tilde{a}_{d-1}+\tilde{c}_{d-1})\A_d.%=\tilde{b}_{d-2}\A_{d-1}+\delta\A_d.
$$
Thus,
$$
\A\tilde{p}_{d-1}(\A)+\delta\tilde{p}_{d-1}(\A)=
\tilde{b}_{d-2}\A_{d-1}+\delta\A_d.
$$
Now, the key point is that the distance matrix $\A_i$ of a bipartite graph can be thought as a $2\times 2$ block matrix such that when $i$ is odd (respectively, even) the diagonal (respectively, off-diagonal) entries are the zero matrices. Indeed, the same happens with the powers of the adjacency matrix, $\A^{\ell}$. Assume that $d$ is odd (respectively, even), the odd part (respectively, even part) of $x\tilde{p}_{d-1}+\delta\tilde{p}_{d-1}$ is a polynomial $p$ satisfying $p(\A)=\delta\A_d$ (the other part gives $\tilde{b}_{d-2}\A_{d-1}$). Hence, from Theorem \ref{teo caract d-r}, $\G$ is distance-regular, with highest degree distance polynomial $p_d=\frac{1}{\delta}p$.
\final

To illustrate the above result, let us give two examples: one for the bipartite case and the other for the nonbipartite one. With respect to the former,
consider the cube $Q_3$, which satisfies the conditions of the theorem and has the following parameters:
\begin{itemize}
\item
$\iota(Q_3)=\{3,2,1;1,2,3\}$,\quad $p_0=1$,\ \ $p_1=x$,\ \ $p_2=\frac{1}{2}(x^2-3)$,\ \ $p_3=\frac{1}{6}(x^3-7x)$;
\item
$\tilde{\iota}(Q_3)=\{2,1;1,2\}$,\qquad\quad $\tilde{p}_0=1$,\ \ $\tilde{p}_1=x-1$,\ \ $\tilde{p}_2=\frac{1}{2}(x^2-2x-1)$.
\end{itemize}
Then, as $d=\delta=3$ and $\tilde{b}_{1}=1$, we have
$$
\textstyle
x\tilde{p}_{2}+\delta\tilde{p}_{1}=\frac{1}{2}(x^3-7x)+\frac{1}{2}(x^2-3)=\frac{1}{3}p_3+p_2,
$$
as it should be.

Now, for the nonbipartite case we consider the odd graph $O_4$, which also satisfies the conditions of the theorem, and has parameters:
\begin{itemize}
\item
$\iota(O_4)=\{4,3,3;1,1,2\}$,\quad $p_0=1$,\ \ $p_1=x$,\ \ $p_2=x^2-4$,\ \ $p_3=\frac{1}{2}(x^3-7x)$;
\item
$\tilde{\iota}(O_4)=\{3,3,2;1,1,4\}$,\quad $\tilde{p}_0=1$,\ \ $\tilde{p}_1=x-1$,\ \ $\tilde{p}_2=x^3-x-3$,
\item[] \hskip 4.3cm $\tilde{p}_3=\frac{1}{4}(x^3-2x^2-5x+6)$.
\end{itemize}
Thus, with $i=2$ and since $d=3$ and $\delta=4$, Eq. (\ref{p_(i+1)}) gives
$$
\textstyle
p_{3}=\frac{1}{\tilde{a}_2+c_2}((x+\delta)\tilde{p}_{2}-\tilde{b}_{1}p_2)
=\frac{1}{2}((x+4)(x^2-x-3)-3(x^2-4))= \frac{1}{2}(x^3-7x),
$$
as claimed.

Using the algebraic approach, we can also derive some properties of the edge-intersection numbers,
apart from the trivial one $\tilde{a}_i+\tilde{b}_i+\tilde{c}_i=\delta$. As an example, we have the following:

\begin{lema}
Let $\G$ be a nonbipartite edge-distance-regular graph with diameter $d$. Then,
\begin{equation}\label{relation-ein}
\tilde{a}_{i}=\tilde{b}_{i-1}-\tilde{b}_{i},\qquad i=1,\ldots,d-1.
\end{equation}
\end{lema}
\prova
In the proof of Theorem \ref{thm edr is dr algebraic}, we have already shown that, if $w\in V_{i,i-1}$ then, $(\B_{i}\B_0^{\top})_{uw}=\tilde{b}_{i-1}$. Similarly, one can prove that, if $w\in V_{i,i+1}$,
then $(\B_{i}\B_0^{\top})_{uw}=\tilde{a}_{i}+\tilde{b}_{i}$. Since, in both cases, $\dist(u,w)=i$,
the two values must be equal and the result follows. \final

%\section{Problems}
%\begin{enumerate}
%\item
%Prove or give a counterexample: Every {\em edge-distance-regularised} (that, is, edge-distance-regular around each of its edges) is edge-distance-regular.
%(See \cite[Th. 4.7]{cdffg11} for an spectral excess theorem characterizing edge-distance-regularised graphs.)
%\item
%Prove Theorem \ref{thm edr is dr} by using the corresponding versions of the spectral excess theorem. (To do so, we need to relate the local spectra of a vertex and an edge.)
%\item
%Idem by using the characterizations based on the invariance of the numbers of walks of a given length (between vertices---for drg's---or between an edge and a vertex---for edrg's).
%\item
%Since we know that every bipartite distance-regular graph is edge-distance-regular, could we use the theory of the later to provide new insight into the former? (For instance, could we define some Klein-like parameters for edge-distance-regularity in order to have extra conditions for feasibility of possible arrays?)
%
%\end{enumerate}

\vskip 1cm
\noindent{\large \bf Acknowledgments.}  Research supported by the
{\em Ministerio de Ciencia e Innovaci\'on}, Spain, and the {\em European Regional
Development Fund} under project MTM2011-28800-C02-01$^b$, the {\em Catalan Research
Council} under project 2009SGR1387$^b$, and by the Netherlands Organization of Scientific Research (NWO)$^a$.
The authors are most grateful to Prof. E.~Garriga for his helpful comments.

\newpage
%%%%%%%%%%%%%%%%%%%%%%%%%%%%%%%%%%%%%%%%%%%%%%%%
%Bibliografia
%%%%%%%%%%%%%%%%%%%%%%%%%%%%%%%%%%%%%%%%%%%%%%%%


\begin{thebibliography}{99}

\bibitem{b93}
N. Biggs, \emph{Algebraic Graph Theory}, Cambridge University Press,
Cambridge, 1974, second edition, 1993.

\bibitem{bcn89}
A.E. Brouwer, A.M. Cohen, and A. Neumaier, \emph{Distance-Regular Graphs},
Springer-Verlag, Berlin-New York, 1989.

%\bibitem{cffg09}
%M. C\'amara, J. F\`abrega, M.A. Fiol, and E. Garriga,
%Some families of orthogonal polynomials of a discrete variable and
%their applications to graphs and codes, {\em Electron. J. Combin.} {\bf 16(1)} (2009), \#R83.

\bibitem{cdffg11}
M. C\'amara, C. Dalf\'o, J. F\`abrega, M.A. Fiol, and E. Garriga, Edge-distance-regular
graphs, {\em J. Combin. Theory Ser. A} {\bf 118} (2011), 2071--2091.

\bibitem{cds80}
C.D. Cvetkovi\'c, M. Doob, and H. Sachs, {\it Spectra of Graphs}, third edition, Johann Barth Verlag, 1995. First edition: Deutscher Verlag der Wissenschaften, Berlin, 1980, Academic Press, New York, 1980.

\bibitem{ddfgg11}
C. Dalf\'{o}, E.R. van Dam, M.A. Fiol, E. Garriga, and B.L.
Gorissen, On almost distance-regular graphs, {\em J. Combin.
Theory Ser. A} \textbf{118} (2011), 1094--1113.

\bibitem{cn05}
B. Curtin, and K. Nomura, 1-homogeneous, pseudo-1-homogeneous, and 1-thin distance-regular graphs. {\it J. Combin. Theory Ser. B} {\bf 93} (2005), 279--302.

\bibitem{vd08}
E.R. van Dam, The spectral excess theorem for distance-regular
graphs: a global (over)view, {\em Electron. J. Combin.} {\bf 15(1)}
(2008), \#R129.

\bibitem{vdf12}
E.R. van Dam, and M.A. Fiol,
A short proof of the odd-girth theorem, {\em Electron. J. Combin.} {\bf 19(3)}
(2012), \#P12.

\bibitem{vdh10}
E.R. van Dam, and W.H. Haemers, An odd characterization of the generalized odd graphs, {\em J. Combin. Theory, Ser. B} {\bf 101} (2011), no. 6, 486--489.

\bibitem{f02}
M.A. Fiol, Algebraic characterizations of
distance-regular graphs, {\it Discrete Math.} {\bf 246} (2002), 111--129.

\bibitem{fgg09}
M.A. Fiol, S. Gago, and E. Garriga, A simple proof of the spectral excess theorem for distance-regular graphs, {\it Linear Algebra Appl.} {\bf 432} (2010), 2418--2422.

\bibitem{fg97}
M.A. Fiol, and E. Garriga, From local adjacency
polynomials to locally pseudo-distance-regular graphs, {\it
J.  Combin. Theory Ser. B} {\bf 71} (1997), 162--183.

\bibitem{fg01}
M.A. Fiol, and E. Garriga, An algebraic characterization of completely regular codes in distance-regular
graphs, {\it SIAM J. Discrete Math.} {\bf 15} (2001), no. 1, 1--13.

\bibitem{fgy96b}
M.A. Fiol,  E. Garriga, and J.L.A. Yebra, Locally
pseudo-distance-regular graphs, {\it J.  Combin. Theory Ser. B}
{\bf 68} (1996), 179--205.

\bibitem{g93}
C.D. Godsil, {\it Algebraic Combinatorics}, Chapman and Hall, New York, 1993.

\bibitem{gmk80}
C.D. Godsil, and B.D. McKay,
Feasibility conditions for the existence of walk-regular graphs,
{\em Linear Algebra Appl.} {\bf 30} (1980), 51--61.

\bibitem{hof63}
A.J. Hoffman, On the polynomial of a graph, {\it Amer. Math. Monthly} {\bf 70} (1963), 30--36.

\bibitem{lw11}
G.-S. Lee and C.-w. Weng, The spectral excess theorem for
general graphs, {\em J. Combin. Theory, Ser. A} 119 (2012), 1427--1431.

\bibitem{m92}
W.J. Martin, {\em Completely Regular Subsets}, Ph.D. Thesis,
University of Waterloo, Waterloo, 1992. (Available in {\tt http://users.wpi.edu/$\thicksim$martin/RESEARCH/THESIS/}.)

\bibitem{n94}
K. Nomura, Homogeneous graphs and regular near polygons, {\it J.  Combin. Theory Ser. B} {\bf 60} (1994), 63--71.


\end{thebibliography}
\end{document}